\documentstyle[12pt]{article}

\def\hybrid{\topmargin 0pt      \oddsidemargin 0pt
        \headheight 0pt \headsep 0pt
        \textwidth 6.35in       
        \textheight 9.5in       
        \marginparwidth 0.0in
        \parskip 5pt plus 1pt   \jot = 1.5ex}
\catcode`\@=11
\def\marginnote#1{}

\newcount\hour
\newcount\minute
\newtoks\amorpm
\hour=\time\divide\hour by60
\minute=\time{\multiply\hour by60 \global\advance\minute by-\hour}
\edef\standardtime{{\ifnum\hour<12 \global\amorpm={am}%
        \else\global\amorpm={pm}\advance\hour by-12 \fi
        \ifnum\hour=0 \hour=12 \fi
        \number\hour:\ifnum\minute<10 0\fi\number\minute\the\amorpm}}
\edef\militarytime{\number\hour:\ifnum\minute<10 0\fi\number\minute}

\def\draftlabel#1{{\@bsphack\if@filesw {\let\thepage\relax
   \xdef\@gtempa{\write\@auxout{\string
      \newlabel{#1}{{\@currentlabel}{\thepage}}}}}\@gtempa
   \if@nobreak \ifvmode\nobreak\fi\fi\fi\@esphack}
        \gdef\@eqnlabel{#1}}
\def\@eqnlabel{}
\def\@vacuum{}
\def\draftmarginnote#1{\marginpar{\raggedright\scriptsize\tt#1}}

\def\draftlabel#1{{\@bsphack\if@filesw {\let\thepage\relax
   \xdef\@gtempa{\write\@auxout{\string
      \newlabel{#1}{{\@currentlabel}{\thepage}}}}}\@gtempa
   \if@nobreak \ifvmode\nobreak\fi\fi\fi\@esphack}
        \gdef\@eqnlabel{#1}}
\def\@eqnlabel{}
\def\@vacuum{}
\def\draftmarginnote#1{\marginpar{\raggedright\scriptsize\tt#1}}

\def\draft{\oddsidemargin -.5truein
        \def\@oddfoot{\sl preliminary draft \hfil
        \rm\thepage\hfil\sl\today\quad\militarytime}
        \let\@evenfoot\@oddfoot \overfullrule 3pt
        \let\label=\draftlabel
        \let\marginnote=\draftmarginnote
   \def\@eqnnum{(\theequation)\rlap{\kern\marginparsep\tt\@eqnlabel}%
\global\let\@eqnlabel\@vacuum}  }

\def\numberbysection{\@addtoreset{equation}{section}
        \def\theequation{\thesection.\arabic{equation}}}

\def\underline#1{\relax\ifmmode\@@underline#1\else
        $\@@underline{\hbox{#1}}$\relax\fi}

\def\titlepage{\@restonecolfalse\if@twocolumn\@restonecoltrue\onecolumn
     \else \newpage \fi \thispagestyle{empty}\c@page\z@
        \def\thefootnote{\fnsymbol{footnote}} }

\def\endtitlepage{\if@restonecol\twocolumn \else  \fi
        \def\thefootnote{\arabic{footnote}}
        \setcounter{footnote}{0}}  
\relax

\numberbysection
\hybrid

\def\beq{\begin{equation}}
\def\eeq{\end{equation}}
\def\p{\partial}

\newtheorem{th}{Theorem}[section]
\newtheorem{prop}{Proposition}[section]
\newtheorem{cor}{Corollary}[section]
\newtheorem{lem}{Lemma}[section]

\def\square{\hfill
{\vrule height6pt width6pt depth1pt} \break \vspace{.01cm}}

\begin{document}

\begin{titlepage}

\title{On the spectral curve of the difference Lam\'e operator}

\author{A. Zabrodin
\thanks{Joint Institute of Chemical Physics, Kosygina str. 4, 117334,
Moscow, Russia and ITEP, 117259, Moscow, Russia}}
\date{December 1998}

\maketitle

\begin{abstract}

We give two "complementary" descriptions of the curve $\Gamma$
parametrizing double-\-Bloch solutions to the difference
analogue of the Lam\'e equation. The curve depends on
a positive integer number $\ell$ and two continuous
parameters: the "lattice spacing" $\eta$ and the modular
parameter $\tau$. Apart from being a covering of the
elliptic curve with the modular parameter $\tau$,
$\Gamma$ is a hyperelliptic curve of genus $2\ell$.
We also point out connections between the spectral curve
and representations of the Sklyanin algebra.

\end{abstract}


\end{titlepage}

\section{Introduction}

Schr\"odinger operators with a periodic potential
usually have infinitely many gaps in the spectrum.
Exceptional cases, when there are only a finite number
of gaps, are of particular interest for the theory of ordinary
differential equations as well as for applications. Their
study goes back to classical works of the last century.
The renewed interest to the theory of finite-gap operators
is due to their role
in constructing quasi-periodic exact solutions
to non-linear integrable equations.

Among known examples of the finite-gap operators, the most
familiar one is the classical Lam\'e operator
\beq
\label{Int0}
{\cal L}=-{d^2\over dx^2}+\ell (\ell +1)\wp(x)\,,
\eeq
where $\wp(x)$ is the Weierstrass $\wp$-function and
$\ell$ is a parameter. The finite gap property of
higher Lam\'e operators for integer values of
$\ell$ was established in \cite{ince}. If $\ell$ is a positive
integer, then there exists a differential operator of order
$2\ell +1$ that commutes with ${\cal L}$, so
the Lam\'e operator has exactly
$\ell$ gaps in the spectrum.
Such a remarkable spectral property is a
signification of a hidden algebraic
symmetry underlying the spectral problem.

In \cite{kz}, a connection
between the finite-gap integration theory of soliton equations and
the representation theory of Sklyanin algebra \cite{Skl1}
was found and the following
difference analogue of the Lam\'e operator was proposed:
\beq
L = \frac{\theta _{1}(2x -2\ell \eta)}
{\theta _{1}(2x)}\, e^{\eta \p _x}
+\frac{\theta _{1}(2x +2\ell \eta)}
{\theta _{1}(2x)}\, e^{-\eta \p _x}\,.
\label{Int1}
\eeq
Here $\theta_1(x)\equiv \theta_1(x|\tau)$ is the odd
Jacobi $\theta$-function, $\ell$ is a non-negative integer
and $\eta \in {\bf C}$ is a parameter which is assumed to
belong to the fundamental parallelogram
with vertices $0$, $1$, $\tau$, $1+\tau$.
The origin of the operator (\ref{Int1}) is traced back
to Sklyanin's paper \cite{Skl2} of 1983, where a functional
realization of the Sklyanin algebra was found. Namely,
$L$ coincides with one
of the four generators of the Sklyanin algebra in the
functional realization. Therefore, the Sklyanin algebra
provides a natural algebraic framework for analyzing the
spectral properties of the operator $L$.
(A different algebraic approach to the difference
analogues of the Lam\'e operators was proposed in \cite{FV1}.)
Nowdays, another face of this operator is probably more
familiar: it is the hamiltonian of the elliptic
two-body Ruijsenaars model \cite{Ruij}.

As is already expected from the relation to the Sklyanin
algebra, the spectral problem $L\Psi =E\Psi$
is closely connected
with the simplest one-site $XYZ$ spin chain of spin $\ell$ at
the site. Indeed, the operator $L$ is proportional to
the trace of the quantum $\mbox{{\sf L}}$-operator
of this model.
Integrable spin chains of $XYZ$-type can be solved
by the generalized algebraic Bethe
ansatz \cite{Baxter},\cite{FT},\cite{Takebe}.
In our case the Bethe ansatz approach amounts to
looking for the eigenfunctions of the form
$$
\Psi (x)=K^{x/\eta}\prod_{j=1}^{\ell}\theta_1(2x-2x_j)\,,
$$
where $K$ and $x_j$ are parameters. If they are constarined
by the system of Bethe equations
$$
K^2 \frac{\theta_1(2x_i-2\ell \eta)}{\theta_1(2x_i+2\ell \eta)}
=\prod _{j=1, \neq i}^{\ell}
\frac{\theta_1(2x_i-2x_j -2\eta)}{\theta_1
(2x_i-2x_j +2\eta)}\,,
\;\;\;\;\;\; i=1,\ldots , \ell\,,
$$
then $\Psi$ is an eigenfunction of $L$.
The energy $E$ obtained from the eigenvalue equation at
a particular value of $x$ (say, $x=\ell \eta$),
$$
E=K^{-1}\frac{\theta_1(4\ell \eta)}{\theta_1(2\ell \eta)}
\prod_{j=1}^{\ell}\frac{\theta_1(2(\ell -1)
\eta -2x_j)}{\theta_1(2\ell \eta -2x_j)}\,,
$$
is then a multivalued function of $K$. It becomes single-valued
on the {\it spectral curve} $\Gamma$ of the operator $L$,
which, therefore, carries all the information about its spectral
properties. The points $P$ of the curve are solutions
$P=\{K,x_1, \ldots , x_{\ell}\}$ to the Bethe system.
However, the description of the curve provided by
the Bethe equations is neither the
most economic nor very informative one (at least for small
values of $\ell$).

Let $\Psi (x)$ be a function on which the operator
(\ref{Int1}) acts.
Putting $\Psi_n=\Psi (n\eta+x_0)$,
we assign to (\ref{Int1}) the difference Schr\"odinger operator
$L \Psi_n=A_n \Psi_{n+1}+B_n \Psi_{n-1}$
with quasiperiodic coefficients. The spectrum of a generic
operator of this form
has a structure of Cantor set type. If $\eta$ is a rational
number, $\eta=P/(2Q)$, this operator
has $Q$-periodic coefficients.
In general, $Q$-periodic difference Schr\"odinger
operators have $Q$ stable bands in the spectrum.

It was proved \cite{kz}
that for integer $\ell$ the operator $L$ is
algebraically integrable and, therefore,
is a difference analogue of the classical Lam\'e
operator which can be obtained
from $L$ in the limit $\eta\to 0$.
Algebraic integrability of $L$ implies, in particular,
some extremely unusual spectral properties of this operator.
Namely, the operator $L$ given by eq.\,(\ref{Int1})
for positive integer values
of $\ell$ and {\it arbitrary} generic $\eta$
has $2\ell +1$ stable bands (and $2\ell$ gaps) in the spectrum.
Its Bloch functions are parametrized by points $P=(w,E)$
of a hyperelliptic curve of genus $2\ell$ defined by the equation
\beq
w^2=c \prod_{i=1}^{2\ell+1}(E^2-E_i^2)\,,
\label{15a}
\eeq
where $c$ is a constant (introduced here for consistency
with the definition of $w$ given below).
Moreover, eigenfunctions of the operator $L$ at the edges
of bands $E=\pm E_{i}$ span an invariant functional
subspace for all generators of the Sklyanin algebra.

The finite-gap property of the operator $L$ means that
there exists a difference operator $W$ of finite order
such that it can not be rerpesented as a polynomial
function of $L$, and that commutes with $L$: $[L,W]=0$.
(This is the
difference version of the Novikov equation
for coefficients of the operators $L$ and $W$.)
It is well-known from the early days of
finite-gap theory that the ring of operators commuting with the
finite-gap operator is isomorphic to a ring of meromorphic functions
on the corresponding spectral curve with poles at "infinite
points". For difference operators
this was proved in
\cite{mum},\,\cite{kr1}. Therefore, the ring of
operators commuting with $L$
is generated by $L$ and an operator $W$ such
that
\beq
\label{15b}
W^2=c\prod_{i=1}^{2\ell +1}(L^{2}-E_i^2)\,.
\eeq
The variable $w$ in (\ref{15a}) is eigenvalue of
the operator $W$, i.e., $W\Psi =w\Psi$, where
$\Psi$ is a common eigenfunction of $L$ and $W$.
The hyperelliptic curve (\ref{15a}) has two
"infinite points" $\infty _{\pm}$, where the function
$E$ has first order poles.
In \cite{krnov}, for any algebraic curve
with two punctures,
a special basis in the ring
of meromorphic functions with
poles at the punctures was introduced.
Due to the isomorphism between the ring of meromorphic
functions with poles at $\infty _{\pm}$ and
the ring of commuting operators, there exist
operators commuting with $L$ such that their
eigenvalues on the common (Baker-Akhiezer)
eigenfunction coincide with the basis functions.
The explicit form of
these operators and the operator
$W$ in the case when $L$ is
the difference Lam\'e operator (\ref{Int1}) was
found in \cite{FV2}. Remarkably, this commuting
family (parametrized by a complex parameter) coincides
with the Baxter $Q$-operator for the one-site $XYZ$-model
with spin $\ell$.

This paper is devoted to a more detailed analysis of
the spectral curves of the difefference Lam\'e operators
for arbitrary positive integer values of $\ell$.
Let us mention that spectral curves of the
classical Lam\'e operator (\ref{Int0}) and its
Treibich-Verdier generalizations \cite{TV} for small
values of $\ell$ were studied in \cite{En}.
Section 2 contains some algebraic preliminaries
on quantum transfer matrices and the $Q$-operator
for the simplest one-site $XYZ$-model with spin
$\ell \in {\bf Z}_{+}$.
This algebraic framework is very helpful since it
allows one to represent the equation of the spectral curve
in the most compact explicit form. The exposition in
the first part of Section 3 follows that of the paper \cite{kz}.
Namely, we construct double-Bloch solutions to the difference
Lam\'e equation and obtain the spectral curve defined by
two equations for three variables. One of these variables
is the eigenvalue $E$, the other two parametrize the Bloch
multipliers of the solution. Next, we show that $E$ can
be eliminated leaving us with one equation for two variables.
In this realization, the fact that the spectral curve
is a covering of the elliptic curve is transparent.
However, the hyperelliptic property of the curve is implicit.
In Section 4 we present the curve
in the explicit hyperelliptic form. Finally, Section 5 contains
a few results on the connection with representations
of the Sklyanin algebra.

\section{Quantum transfer matrices
and the $Q$-operator for the one-site $XYZ$-model with spin
$\ell$}

This section contains selected ingredients
of the quantum inverse
scattering approach to $XYZ$ spin chains \cite{FT}
specified to the case of the one-site "chain"
with spin $\ell$. These constructions
turn out to be particularly
useful for analyzing the spectral
curve of the difference Lam\'e operator.

We begin with a few formulas related to
the Sklyanin algebra and its representations.
Definitions and transformation properties
of the Jacobi $\theta$-functions
$\theta_a(x|\tau )$ are listed in Appendix A.
For brevity, we write $\theta _a(x|\tau )
\equiv \theta _a(x)$.

The elliptic quantum $\mbox{{\sf L}}$-operator is
the matrix
\beq
\mbox{{\sf L}}(u)=\frac{1}{2}
\left ( \begin{array}{cc}
\theta _{1}(u){\cal S}_0 +\theta _{4}(u){\cal S}_3 &
\theta _{2}(u){\cal S}_1 +\theta _{3}(u){\cal S}_2
\\& \\
\theta _{2}(u){\cal S}_1 -\theta _{3}(u){\cal S}_2 &
\theta _{1}(u){\cal S}_0 -\theta _{4}(u){\cal S}_3
\end{array} \right )
\label{L}
\eeq
with non-commutative matrix elements. Specifically,
${\cal S}_a$ are difference operators in a
complex variable $x$:
\beq
{\cal S}_{a} = \frac{\theta _{a+1}(2x -2\ell \eta)}
{\theta _{1}(2x)}\, e^{\eta \p _x}
-\frac{\theta _{a+1}(-2x -2\ell \eta)}
{\theta _{1}(2x)}\, e^{-\eta \p _x}\,,
\label{Sa}
\eeq
introduced by Sklyanin \cite{Skl2} in 1983.
Comparing with (\ref{Int1}), we identify $L={\cal S}_0$.
The four operators ${\cal S}_a$ obey the commutation
relations of the Sklyanin algebra\footnote{The standard generators
of the Sklyanin algebra \cite{Skl1}
are related to ours as follows:
$S_{a}=(i)^{\delta _{a,2}}\theta _{a+1}(\eta ){\cal S}_a$.}:
\beq
\begin{array}{l}
(-1)^{\alpha +1}I_{\alpha 0}{\cal S}_{\alpha}{\cal S}_{0}=
I_{\beta \gamma}{\cal S}_{\beta}{\cal S}_{\gamma}
-I_{\gamma \beta}{\cal S}_{\gamma}{\cal S}_{\beta}\,,
\\ \\
(-1)^{\alpha +1}I_{\alpha 0}{\cal S}_0 {\cal S}_{\alpha}=
I_{\gamma \beta}{\cal S}_{\beta}{\cal S}_{\gamma}
-I_{\beta \gamma}{\cal S}_{\gamma}{\cal S}_{\beta}
\end{array}
\label{skl6}
\eeq
with the structure constants
$I_{ab}=\theta _{a+1}(0)\theta _{b+1}(2\eta)$.
Here $a,b =0, \ldots , 3$ and
$\{\alpha ,\beta , \gamma \}$ stands for any cyclic
permutation of $\{1 ,2,3\}$.
The relations of the Sklyanin algebra
are equivalent to the condition that the $\mbox{{\sf L}}$-operator
satisfies the
"$\mbox{{\sf R}}\mbox{{\sf L}}
\mbox{{\sf L}}=\mbox{{\sf L}}
\mbox{{\sf L}}\mbox{{\sf R}}$"
relation with the
standard elliptic $R$-matrix.

The parameter
$\ell$ in (\ref{Sa}) is called spin of the representation.
If necessary, we write ${\cal S}_a ={\cal S}_{a}^{(\ell )}$
to indicate the dependence on $\ell$.
When $\ell \in
\frac{1}{2} {\bf Z}_{+}$, these operators
have a finite-dimensional invariant subspace, namely,
the space
$\Theta_{4\ell}^{+}$ of {\it even} $\theta$-functions of
order $4\ell$ (see Appendix A). This
is the representation space of the $(2\ell +1)$-dimensional
irreducible representation (of series a))
of the Sklyanin algebra.

Trace of $\mbox{{\sf L}}(u)$
(in the two-dimensional auxiliary space),
that is the simplest
quantum transfer matrix $T_1(u)$, is proportional
to $L={\cal S}_0$:
\beq
\label{T1}
T_1(u)=\mbox{tr}\mbox{{\sf L}}(u)=\theta _{1}(u){\cal S}_0\,.
\eeq
The whole family of commuting transfer matrices
$T_s(u)$, $s\in {\bf Z}_{+}$, is
obtained from (\ref{L}) via the fusion procedure \cite{KRS}.
We denote them by $T_s(u)$, $s\in {\bf Z}_{+}$.
They commute for all values of $s$ and $u$:
$[T_s(u),\, T_{s'}(u')]=0$.
Here we do not need to recall the fusion procedure itself
and refer the reader to \cite{KR} and \cite{Tak}, where
integrable magnets of higher spin in the $XXZ$ and $XYZ$ case,
respectively, were constructed by means
of the fusion procedure. For our purpose it is enough
to {\it define} the operators $T_s(u)$ by the
recurrence relations (known as the fusion relations):
\beq
\begin{array}{lll}
s<2\ell: &&
\!\!\!\!T_1 (u\!-\!s\eta )T_s(u\!+\!\eta )\!=\!T_{s\!+\!1}(u)\!+\!
\theta _{1}\bigl (u\!-\!s\eta \!-\!2\ell \eta \bigr)
\theta _{1}\bigl (u\!-\!s\eta \!+\!2(\ell \!+\!1)\eta \bigr)
T_{s\!-\!1}(u\!+\!2\eta ), \\ && \\
s=2\ell: &&
\!\!T_1 (u\!-\!2\ell \eta )T_{2\ell}(u\!+\!\eta )\!
=\!\theta_1(u)T_{2\ell \!+\!1}(u)\!+\!
\theta _{1}\bigl (u\!+\!2\eta \bigr)
\theta _{1}\bigl (u\!-\!4\ell \eta \bigr)
T_{2\ell \!-\!1}(u\!+\!2\eta )\,, \\ && \\
s>2\ell: &&
\!\!T_1 (u\!-\!s\eta )T_s(u\!+\!\eta )\!=
\!\theta_1 \bigl (u\!-\!s\eta \!+\!2\ell \eta \bigr )
T_{s\!+\!1}(u)\!+\!
\theta _{1}\bigl (u\!-\!s\eta \!-\!2\ell \eta \bigr)
T_{s\!-\!1}(u\!+\!2\eta )
\end{array}
\label{fus1}
\eeq
with the "initial condition"
$T_0(u)=\mbox{id}$, $T_1(u)=\theta _{1}(u){\cal S}_0$.
There is a useful determinant formula which
solves the recurrence relations and represents
$T_s(u)$ through $T_1(u)$ \cite{BR}:
\beq
\begin{array}{llll}
T_s(u)&\!\!=&\!\!\det ({\bf T}_{ij}(s,u))_{1\leq i,j\leq s}\,,
\;\;\;\;\;\; 0\leq s\leq 2\ell\,, &\\ &&&\\
T_s(u)\!\!&=& \!\!
\left(\displaystyle{\prod_{i=1}^{s-2\ell}}
\theta_1(u\!+\!(2\ell \!+\!2i \!-\!s\!-\!1)\eta )\right )^{-1}
\!\!\det ({\bf T}_{ij}(s,u))_{1\leq i,j\leq s}\,,&
\;\;\;\;\; s> 2\ell\,,
\end{array}
\label{det}
\eeq
where
\beq
\label{Tij}
\begin{array}{lll}
{\bf T}_{ij}(s,u)&=&\delta _{i,j-1}
\theta _{1}\bigl (u+(s\!-\!2\ell \!-\!1\!-\!2i)\eta \bigr )
+\delta _{i,j+1}
\theta _{1}\bigl (u+(s\!+\!2\ell \!+\!3\!-\!2i)\eta \bigr )
\\ && \\
&+&\delta _{i,j}T_1\bigl (u +(s\!+\!1\!-\!2i)\eta \bigr )\,.
\end{array}
\eeq

Let $\Psi$ be a common eigenfunction of $L={\cal S}_0$
and $T_s(u)$, and let $E$ be the eigenvalue of $L$:
$L\Psi =E\Psi$. Then the eigenvalue of $T_s(u)$ is
a polynomial in $E$ of degree $s$, which we denote by
$T_s(u,E)$:
\beq
\label{TE}
\left \{ \begin{array}{lll}
{\cal S}_0 \Psi &=&E\Psi \,,\\
T_s(u)\Psi &=& T_s(u,E)\Psi \,.
\end{array} \right.
\eeq
The eigenvalues $T_s(u,E)$ are determined by eqs.\,(\ref{det}),
(\ref{Tij}), where $T_1(u\!+\!(s\!+\!1\!-\!2i)\eta)$ is replaced
by $E\theta_1(u\!+\!(s\!+\!1\!-\!2i)\eta)$.
As a function of $u$, $T_s(u,E)$ for
$1\leq s \leq 2\ell$ is easily seen to belong
to the space $\Theta_s$ of $\theta$-functions of order $s$
(for the precise definition see Appendix A).
An important fact,
not obvious from the definition (\ref{fus1}),
is that $T_s(u,E)$ for all
$s>2\ell$ belong to the space $\Theta_{2\ell}$, i.e.,
the denominator in (\ref{det}) cancels. In particular,
\beq
\label{T2ell1}
\begin{array}{lll}
T_{2\ell +1}(u,E)&=&\displaystyle{\frac{1}{\theta _{1}(u)}}
\det \Bigl (\delta_{i,j}E\theta_1 (u+2(\ell +1-i)\eta )+ \\
&&\;\;\;\;+\,\delta_{i,j-1}\theta_1 (u\!+\!2(2\ell
\!+\!1\!-\!i)\eta )+
\delta_{i,j+1}\theta_1 (u\!-\!2(i\!-\!1)\eta )
\Bigr )_{1\leq i,j \leq 2\ell +1}
\end{array}
\eeq
is holomorphic at $u=0$.

The full family of operators commuting with ${\cal S}_0$
is generated
by Baxter's $Q$-operator $\hat Q(u)$. Moreover, the operators
$\hat Q(u)$ commute with all the transfer matrices and
among themselves: $[T_s(u),\, \hat Q(u')]=0$, $[\hat Q(u),
\, \hat Q(u')]=0$.
They obey the famous Baxter $T$-$Q$-relation
\beq
\label{bax}
\theta _{1}(u-2\ell \eta )\hat Q(u+2\eta )
+\theta _{1}(u+2\ell \eta )\hat Q(u-2\eta )
=T_1(u) \hat Q(u)\,.
\eeq
We also recall
the formula for $T_s(u,E)$ through
eigenvalues $Q(u)$ of the
$Q$-operator \cite{KR},\,\cite{Tak}:
\begin{eqnarray}
\label{00}
T_s(u,E)&=& \frac{Q(u+(s+1)\eta )Q(u-(s+1)\eta )}
{\prod _{p=1}^{2\ell -s}\theta _{1}\Bigl (
u+(2\ell +1 -s-2p)\eta \Bigr )}
\nonumber\\
&\times& \sum _{j=0}^{s}\frac{\prod _{q=1}^{2\ell }
\theta _{1}\Bigl (
u+(2\ell +1 +s-2j-2q)\eta \Bigr )}
{Q(u+(s -2j -1)\eta )Q(u+(s -2j +1)\eta )}\,,
\;\;\;\;\;1\leq s \leq 2\ell\,.
\end{eqnarray}
If $s\geq 2\ell$, there is no denominator in the prefactor.

Let
$\ell \in {\bf Z}_{+}$. In this case
a commuting family, equivalent to the $Q$-operator,
was explicitly constructed
in \cite{FV2}. Consider the operators
\beq
\label{op1}
A_{\lambda}=\sum _{k=0}^{\ell}
A_{k}(x,\lambda)
\,e^{(2k\eta -\ell \eta +\lambda)\p _x}\,,
\eeq
where $\lambda \in {\bf C}$ and
\begin{eqnarray}
\label{op2}
\!\!\!A_{k}(x,\lambda)\!\!
&\,=\,(-1)^k
\displaystyle{\frac{[\ell ]!}{[2\ell ]!}}
\left [\begin{array}{c}\ell \\ k\end{array} \right ]
&\,\prod _{j=0}^{\ell -k-1}
\frac{\theta _{1}\bigl (2x +2(\ell \!-\!j)\eta \bigr )
\theta _{1}\bigl (2\lambda +2(\ell \!-\!j)\eta \bigr )}
{\theta _{1}\bigl (2x +2\lambda +2(k \!-\!j)\eta \bigr )}
\nonumber\\
&&\!\!\!\times\,\,\prod _{j=0}^{k-1}
\frac{\theta _{1}\bigl (2x -2(\ell \!-\!j)\eta \bigr )
\theta _{1}\bigl (2\lambda -2(\ell \!-\!j)\eta \bigr )}
{\theta _{1}\bigl (2x +2\lambda +2(k \!+\!j\!-\!\ell )
\eta \bigr )}\,.
\end{eqnarray}
If $k=0$ or $k=\ell$, the second (respectively, the first)
product is absent. Here and below we use the
"elliptic factorial" and "elliptic binomial" notation:
\beq
\label{binom}
\begin{array}{l}
\displaystyle{[n]!=\prod_{j=1}^{n}[j]}\,,
\;\;\;\;\;\;[j]\equiv \theta_1(2j\eta)\,,
\\ \\
\left [ \begin{array}{c}n\\m\end{array}\right ]
\equiv \displaystyle{\frac{[n]!}{[m]![n-m]!}}\,.
\end{array}
\eeq

\noindent
The main property
of the operators (\ref{op1}) proved in \cite{FV2} is
their commutativity for all values of $\lambda$:
$[A_{\lambda},\,A_{\lambda'}]=0$.
For generic $\lambda$ the chain of shifts in (\ref{op1})
starts from $-\ell \eta +\lambda$ and
all the $\ell \!+\!1$ coefficients in (\ref{op1}) are non-zero.
However, for $\lambda =l\eta$, $l =\ell , \ell-1, \ldots ,
-\ell $ only $\ell -|l|+1$ of them are non-zero. (For example,
$A_{\ell \eta}=1$,
$A_{(\ell -1)\eta}=([\ell ]/[2\ell ]){\cal S}_0$.)

By difference operator in the next theorem we mean a finite
sum $\sum_k f_k(x)e^{k\eta \p _x}$ with integer $k$.
(So $A_{\lambda}$
are difference operators, in this sense, only
if $\lambda = m\eta$ with integer $m$.)
\begin{th}\cite{FV2}
The ring of difference operators
commuting with
$L={\cal S}_0$ (\ref{Int1}) is
generated by $L$ and $A\equiv A_{(\ell +1)\eta}$.
\end{th}

\noindent
It is convenient to introduce the following special
notation: $A\equiv A_{(\ell +1)\eta}$,
$\bar A\equiv A_{-(\ell +1)\eta}$, $W=A-\bar A$.
The role of the operator $W$ (the very one entering
eq.\,(\ref{15b})) for representations of the Sklyanin
algebra is clarified in Sect.\,5.

We conclude this section by listing some properties
of the operators $A_{\lambda}$ which will be useful in the
sequel.

{\it a) The Baxter $T$-$Q$-relation.}
It has been proved in \cite{FV2} that the $A_{\lambda}$
obey the following operator identity:
\beq
\label{bax1}
{\cal S}_0 A_{\lambda }=
\frac{\theta _{1}(2\lambda -2\ell \eta )}
{\theta _{1}(2\lambda )} A_{\lambda +\eta}
+\frac{\theta _{1}(2\lambda +2\ell \eta )}
{\theta _{1}(2\lambda )} A_{\lambda -\eta}\,,
\eeq
which allows
us to identify $A_{\lambda}$ with the $Q$-operator:
$A_{\lambda}=\hat Q(2\lambda)$
(cf. (\ref{bax}), (\ref{T1})). In other words,
$A_{\lambda}$ can be regarded as an operator solution
to the Baxter relation. The second operator solution
to the second order equation
(\ref{bax1}) is $A_{-\lambda}$. Their
wronskian $\mbox{Wr}(\lambda)$ is
easily evaluated:
\beq
\label{wr}
\mbox{Wr}(\lambda)= A_{\lambda +\eta}A_{-\lambda}-
A_{\lambda}A_{-(\lambda +\eta )}=
([2\ell ]!)^{-1}\left ( \prod _{j=-\ell +1}^{\ell}
\theta_1(2\lambda +2j\eta )\right ) W\,.
\eeq

{\it b) The symmetry $x \leftrightarrow \lambda$.}
Let $F(x)$ be an arbitrary function. Since
$A_k (x,\lambda )=A_k (\lambda , x)$, it is clear from
(\ref{op1}) that the result of action of
the $A_{\lambda}$ on the $F(x)$ is symmetric under the
interchange of $x$ and $\lambda$, i.e.,
\beq
\label{symm}
A_{\lambda}(x, \p _x)F(x)=A_{x}(\lambda ,
\p _{\lambda})F(\lambda)\,,
\eeq
where $A_{\lambda}(x,\p _x)$
(respectively, $A_{x}(\lambda ,\p _{\lambda})$)
acts on the function of $x$ (respectively, of $\lambda$).

{\it c) Even and odd difference operators.}
Let $\Xi$ be the operator changing the sign
of $x$: $\Xi F(x) =F(-x)$. It is clear from the definition that
$\Xi A_{\lambda}\Xi ^{-1} =A_{-\lambda}$.
We call difference operators $O$ such that
$\Xi O\Xi ^{-1} =O$ (respectively,
$\Xi O\Xi ^{-1} =-O$) {\it even} (respectively, {\it odd})
operators. It can be proved \cite{FV2} that for generic
$\eta$ any even
difference operator commuting with
$L$ is a polynomial in $L$. In particular, $A_{k\eta}+A_{-k\eta}$
for $k\in {\bf Z}$ and $A_{\lambda}A_{-\lambda}$
for arbitrary $\lambda ,\eta \in {\bf C}$ are polynomial functions
of $L$.

{\it d) Relations between the transfer matrices $T_s(u)$ and
the difference operators $A_{s\eta}$, $s\in {\bf Z}$.}
From (\ref{op2}) we immediately conclude that
$A_{-j\eta} =A_{j\eta}$,
$-\ell \leq j \leq \ell$, so they are even operators.
Then it follows from the above
that the operators $A_{j\eta}$ with
integer $-\ell \leq j \leq \ell$ are polynomial functions
of ${\cal S}_0$. So, similarly to (\ref{TE}), we define
polynomials $A_{j\eta}(E)$ to be eigenvalues of the
$A_{j\eta}$ on their common eigenfunction
$\Psi$ such that ${\cal S}_{0}\Psi =E\Psi$.
Comparing the fusion relation (\ref{fus1}) for
$s\leq 2\ell$ with (\ref{bax1}),
we identify
\beq
\label{AT}
A_{(\ell -s)\eta }=
\frac{[2\ell -s]!}{[2\ell ]!}
T_s \Bigl ( 2\ell
\eta \!-\!(s\!-\!1)\eta \Bigr )\,,
\;\;\;\;\;\;\;\;
s=0,1, \ldots , 2\ell\,,
\eeq
whence (\ref{det}) yields the determinant
representation of the polynomials $A_{(\ell -s)\eta}(E)$:
\beq
\label{AE}
A_{(\ell -s)\eta}(E)=\left [\begin{array}{c}\ell \\ s
\end{array}\right ]
\left [\begin{array}{c}2\ell \\ s
\end{array}\right ]^{-1} \det \left (
E\delta _{i,j}+
\frac{[-i]}{[\ell \!+\!1\!-\!i]}\delta _{i,j-1}
+\frac{[2\ell \!+\!2\!-\!i]}{[\ell \!+\!1\!-\!i]}\delta _{i,j+1}
\right )_{1\leq i,j\leq s}
\eeq
(here $0\leq s\leq \ell$).
The Baxter equation (\ref{bax1}) gives the recurrence
relation for these polynomials:
\beq
\label{bax2}
A_{(\ell \!-\!s\!-\!1)\eta}(E)=
\frac{[\ell \!-\!s]}{[2\ell \!-\!s]}EA_{(\ell \!-\!s )\eta}(E)
+\frac{[s]}{[2\ell \!-\!s]}A_{(\ell \!-\!s\!+\!1)\eta}(E)
\eeq
with the initial conditions $A_{\ell \eta}(E)=1$,
$A_{(\ell -1)\eta}(E)=([\ell ]/[2\ell ])E$.
It is clear from (\ref{bax2}) that
\beq
\label{parity}
A_{(\ell -s)\eta}(-E)=
(-1)^s A_{(\ell -s)\eta}(E)\,,
\;\;\;\;\;\;
0\leq s\leq \ell\,.
\eeq
At last, we point out the
relation
\beq
\label{TAA}
T_{2\ell +1}(0)=[2\ell ]! \bigl (A_{(\ell +1)\eta } +
A_{-(\ell +1)\eta }\bigr )
\eeq
which follows e.g. from (\ref{00}).
The difference operator in the right hand side
is even. Its eigenvalue is given by the polynomial
$T_{2\ell+1}(0,E)=\lim _{u\to 0} T_{2\ell +1}(u,E)$
(see (\ref{T2ell1})).

\section{Double-Bloch eigenfunctions and explicit
form of the spectral curve}

Let $\ell$ be a positive integer. For our current
purpose it is
more convenient to pass to the function
\beq
\label{db1}
\psi (x) =\Psi (x)\left (
\prod _{j=1}^{\ell}\theta_1(2x-2j\eta)\right )^{-1}.
\eeq
Then the eigenvalue equation for the $L$ acquires the form
\beq
\label{db2}
\psi(x+\eta)+
\frac{\theta_1(2x+2\ell \eta)\theta_1(2x-2(\ell +1)\eta)}
{\theta_1(2x)\theta_1(2x-2\eta)}\,
\psi(x-\eta)=E\psi(x)
\eeq
which we also call the difference analogue of the
Lam\'e equation.
In this form, the coefficient function is
{\it double-periodic} with periods $\frac{1}{2}$ and
$\frac{\tau}{2}$,
so it is natural to look for solutions in the class of
{\it double-Bloch functions} \cite{kz}, i.e., such that
$\psi (x+\frac{1}{2})=B_1 \psi(x)$,
$\psi (x+\frac{1}{2}\tau )=B_{\tau} \psi(x)$ with some
constants $B_1, B_{\tau}$.

Introduce the function
\beq
\label{Phi}
\Phi (x,\zeta)=
\frac{\theta _{1}(\zeta +x)}{\theta _{1}(x)\theta _{1}(\zeta)}\,.
\eeq
Its monodromy properties in $x$ are
$\Phi(x+1, \zeta)=\Phi (x,\zeta)$,
$\Phi(x+\tau, \zeta)=e^{-2\pi i \zeta}\Phi (x,\zeta)$,
i.e., it is a double-Bloch function. The
double-Bloch ansatz for the $\psi$ is
\beq
\label{db3}
\psi(x)=K^{x/\eta}\sum_{j=1}^{\ell}
s_j(\zeta, K, E) \Phi (2x -2j\eta, \zeta)\,,
\eeq
where $\zeta ,K$ parametrize the Bloch multipliers of the
function $\psi (x)$:
$B_1=K^{\frac{1}{2\eta}}$, $B_{\tau}=
K^{\frac{\tau}{2\eta}}e^{-2\pi i \zeta}$.
The coefficients $s_j$ depend on the indicated parameters
only.

Substituting (\ref{db3}) into
(\ref{db2}) and computing the residues at the points
$x=j\eta,\ j=0,\ldots,\ell$, we get $\ell +1$ linear equations
\beq
\sum_{j=1}^{\ell} M_{ij} s_j=0,
\;\;\;\;\;\; i=0,1,\ldots, \ell \,,
\label{db4}
\eeq
for  $\ell$ unknowns $s_j$.
Matrix elements $M_{ij}$ of this system are:

\begin{eqnarray}
\label{db5}
M_{ij}&=&K\delta_{i,j\!-\!1}-E\delta_{i,j}+K^{-1}
\frac{\theta_1(2(j+\ell +1)\eta)
\theta_1(2(j-\ell )\eta)}{\theta_1(2(j+1)\eta)\theta_1(2j\eta)}
\delta_{i,j\!+\!1}+
\nonumber\\
&+&K^{-1}\frac{\theta_1(\zeta \!-\!2(j\!-\!i\!+\!1)
\eta)}{\theta_1(\zeta)}\,
\frac{\theta_1(2(i+\ell)\eta)
\theta_1(2(i\!-\!\ell \!-\!1)\eta)}{\theta_1(2\eta)
\theta_1(2(j\!-\!i\!+\!1)\eta)}
(\delta_{i,0}-\delta_{i,1})\,.
\end{eqnarray}
Here $i=0,1, \ldots , \ell$,
$j=1,2, \ldots , \ell$.
The overdetermined system (\ref{db4}) has nontrivial solutions if
and only if
rank of the rectangular matrix $M_{ij}$ is less
than $\ell$. By $M^{(0)}$ and $M^{(1)}$ we denote
$\ell \times  \ell$ matrices obtained from $M$ by deleting the rows
with $i=0$ and $i=1$,
respectively. Then the set of parameters
$\zeta,K,E$ for which eq.\,(\ref{db4})
has solutions of the form (\ref{db3})
is determined by the system of two equations:
$\det M^{(0)}=\det M^{(1)}=0$.
So the three parameters are constrained by two equations.
They can be written out in a particularly compact form
in terms of the family of polynomials (\ref{AE}) and
the elliptic "binomial coefficients" (\ref{binom}).
Expanding the determinants
with respect to the first row and taking into account
(\ref{AE}), (\ref{bax1}), we come to the following statement.
\begin{th}
The difference Lam\'e equation (\ref{db2})
has double-Bloch solutions of the form (\ref{db3})
if and only if the spectral parameters $\zeta, K, E$
obey the equations
\beq
\begin{array}{l}
\displaystyle{\sum_{j=0}^{\ell}(-1)^j K^{-j}
\theta_1(\zeta-2j\eta)\left [
\begin{array}{c}\ell\\ j \end{array}\right ]}
A_{j\eta}(E)=0\,, \\ \\
\displaystyle{\sum_{j=0}^{\ell +1}(-1)^j K^{-j}
\theta_1(\zeta-2j\eta)\theta_1(2(j-1)\eta)\left [
\begin{array}{c}\ell \!+\!1\\ j \end{array}\right ]}
A_{(j-1)\eta}(E)=0\,,
\end{array}
\label{db7}
\eeq
where $A_{j\eta}(E)$ are polynomials of $(\ell -|j|)$-th degree
explicitly given by (\ref{AE}).
They coincide with eigenvalues of the commuting operators
$A_{j\eta}$ introduced in (\ref{op1}) on their common eigenfunction
$\Psi$ such that $L\Psi = E\Psi$.
\end{th}

The equations (\ref{db7})
define a Riemann surface $\tilde \Gamma$, which covers the
complex plane. The monodromy properties of the $\theta$-function
(see Appendix A) make it clear that
this surface is invariant under the transformation
\beq
\label{db8}
\zeta \longmapsto  \zeta +\tau\,,
\;\;\;\;\;\;
K \longmapsto Ke^{4\pi i\eta}\,.
\eeq
The factor of the $\tilde \Gamma$ over this transformation
is an algebraic curve $\Gamma$, which is
a ramified covering of the elliptic
curve with periods $1$, $\tau$. It is clear
from (\ref{db7}), (\ref{parity}) that the curve
admits the involution
\beq
\label{db81}
(\zeta , K,E)\longmapsto (\zeta , -K,-E)\,,
\eeq
so the spectrum is symmetric with respect
to the reflection $E\rightarrow -E$.
Another result of \cite{kz}, which is not so easy to see
from (\ref{db7}), is that the curve $\Gamma$ is at the same
time a hyperelliptic curve.
\begin{th}\cite{kz}
The curve $\Gamma$
is a hyperelliptic curve of genus $g=2\ell$. The hyperelliptic
involution is given by
\beq
\label{db9}
(\zeta , \,K, \,E) \; \longmapsto \;
(4N\eta -\zeta, \,K^{-1}, E)\,,
\;\;\;\;\;\;N=\frac{1}{2}\ell (\ell +1)\,.
\eeq
The points $P=(\zeta , K, E) \in \Gamma$ of the curve
parametrize double-Bloch solutions
$\psi (x)=\psi (x,P)$ to eq.\,(\ref{db2}), and the
solution $\psi (x,P)$
corresponding to each point
$P\in \Gamma$ is unique up to a constant multiplier.
\end{th}
To compare with \cite{kz}, we note that
the variables $z,k$ used in that paper are
related to $\zeta , K$ as follows:
$\zeta =z +2N\eta$,
$K=k\left (\displaystyle{\frac{\theta _{1}(z-2\eta )}
{\theta _{1}(z+2\eta )}}\right )^{\frac{1}{2}}$.
In terms of $z,k$ the coefficients in
(\ref{db7}) become elliptic functions of $z$
and the hyperelliptic involution is
$z\to -z$, $k\to k^{-1}$.

The edges of bands $\pm E_i$ in (\ref{15a}) are values
of the function $E=E(P)$ at the fixed points of the
hyperelliptic involution.
As is clear from (\ref{db9}), the fixed points lie
above the points $\zeta = 2N\eta +\omega_a$, where $\omega_a$
are the half-periods: $\omega_1=0$, $\omega_2 =\frac{1}{2}$,
$\omega_3 =\frac{1}{2}(1+\tau )$,
$\omega_4 =\frac{1}{2}\tau $. The corresponding values of $K$
are determined from (\ref{db8}).

\begin{cor}
Let ${\cal E}_a$, $a=1,\ldots , 4$ be the set of common
roots of the polynomial equations
\beq
\begin{array}{l}
\displaystyle{\sum_{j=0}^{\ell}
\theta_a(2(N-j)\eta)\left [
\begin{array}{c}\ell\\ j \end{array}\right ]}
A_{j\eta}(E)=0\,, \\ \\
\displaystyle{\sum_{j=0}^{\ell +1}
\theta_a(2(N-j)\eta)\theta_1(2(j\!-\!1)\eta)\left [
\begin{array}{c}\ell \!+\!1\\ j \end{array}\right ]}
A_{(j-1)\eta}(E)=0\,,
\end{array}
\label{db10}
\eeq
where $N=\frac{1}{2}\ell (\ell +1 )$, and
$\theta_a$ are Jacobi $\theta$-functions.
Then the set of the edges of bands $\pm E_i$
is the union of $\,\bigcup_{a=1}^{4}{\cal E}_a$ and its
image under the reflection $E\rightarrow -E$.
\end{cor}

To obtain a more detailed information from
equations (\ref{db7}), one can try to eliminate
$E$ and obtain a single equation
connecting the two Bloch multipliers of the function
(\ref{db3}) (parametrizing through $\zeta$ and $K$).
However, this is not easy to do directly.
A possible way out relies on the following
simple lemma.
\begin{lem}
Let $\Psi (x)$ be any solution to the equation
\beq
\label{eq1}
\frac{\theta _{1}(2x -2\ell \eta)}{\theta _{1}(2x)}
\Psi (x+\eta ) +
\frac{\theta _{1}(2x +2\ell \eta)}{\theta _{1}(2x)}
\Psi (x-\eta )
=E\Psi (x)
\eeq
in the class of entire functions on the complex plane of
the variable $x$, then
\beq
\label{cond}
\Psi (j\eta )=\Psi (-j\eta )\,,
\;\;\;\;\;\;\;\;\; j =1,2, \ldots , \,\ell\,.
\eeq
\end{lem}
This assertion follows from the specific form of the coefficients
of eq.\,(\ref{eq1}). Indeed, putting $x=0$ in (\ref{eq1}),
we have $\Psi (\eta)=\Psi (-\eta)$. The proof can be
completed by induction. At $x=\pm \ell \eta$ one of the coefficients
in the l.h.s. of (\ref{eq1}) vanishes, so
the chain of relations (\ref{cond}) truncates at $j=\ell$.
\square
\paragraph{Remark} The conditions
(\ref{cond}) resemble the "glueing conditions" for the
Baker-Akhiezer function on rational
curves with double points, where they are
imposed on the $\Psi$ with respect to its spectral parameter.
However, contrary to that case, (\ref{cond}) is
imposed in the $x$-plane.

Remarkably, the conditions (\ref{cond})
and the ansatz
\beq
\label{psi1}
\Psi (x)=K^{x/\eta}
\left (\prod _{j=1}^{\ell}\theta _{1}(2x-2j\eta )\right )
\sum _{m=1}^{\ell}s_{m}(K,\zeta )\Phi (2x -2m\eta , \, \zeta )
\eeq
for $\Psi$ (equivalent to the ansatz (\ref{db3})
for $\psi$) with the same function $\Phi (x,z)$ given by
(\ref{Phi}) allow one
to find the relation between the
Bloch multipliers even without explicit use of
the difference Lam\'e equation (\ref{eq1}).
Plugging (\ref{psi1}) into (\ref{cond}), we obtain $\ell$
equalities (for $m=1,2, \ldots , \ell$):
\beq
\label{eq2}
K^m s_m =(-1)^{\ell}K^{-m}\theta _{1}(4m\eta )
\left ( \prod _{j=1, \neq m}^{\ell}
\frac{\theta _{1}(2(m+j)\eta )}
{\theta _{1}(2(m-j)\eta )}\right )
\sum _{n=1}^{\ell} \Phi \Bigl (-2(m+n)\eta , \, \zeta \Bigr )\, s_n\,.
\eeq
This is a system of linear homogeneous equations for $s_n$.
It has nontrivial solutions if and only if its determinant
is equal to zero, whence we obtain
the equation for $\zeta$ and $K$:
\beq
\label{eq3}
\det \Bigl ( K^{2m}\delta _{mn}+G_{mn}(\zeta )
\Bigr )_{1\leq m,n\leq \ell}=0\,,
\eeq
where
$$
G_{mn}(\zeta )=
(-1)^{\ell +1}[2m]
\left ( \prod _{j=1, \neq m}^{\ell}
\frac{[m+j]}{[m-j]}\right )
\Phi \bigl (-2(m+n)\eta , \, \zeta \bigr )\,.
$$
This equation defines a curve $\Gamma _{e}$, which is
the image of the spectral curve $\Gamma$ under the
projection $\Gamma \rightarrow \Gamma _{e}$ that takes
$(\zeta , K, E)$ to $(\zeta , K)$.
A more explicit description of the curve
$\Gamma _{e}$ is given by the following
proposition.
\begin{prop}
The equation of the spectral curve (\ref{eq3})
can be represented in the form
\beq
\label{eq4}
\sum _{j=0}^{N}(-1)^j C^{(\ell)}_{j}
(\eta )\theta _{1}(\zeta -4j\eta)
K^{2(N-j)}=0\,,
\eeq
where $N=\frac{1}{2}\ell (\ell +1)$ and $C_{j}^{(\ell)}(\eta)$ are
coefficients depending only on $\eta$ (and $\tau$) such that
$C_{j}^{(\ell)}(\eta)=C_{N-j}^{(\ell)}(\eta)$,
$C_{0}^{(\ell)}(\eta)=1$.
\end{prop}
{\it Proof.} The following direct proof allows us to find
the explicit form of the coefficients $C_{j}^{(\ell)}$.
To expand the determinant in powers of $K$, we make use
of the identity
$$
\det \left (\frac{\theta_{1}
(x_i+x_j+\zeta)}{\theta_{1}(x_i+x_j)}
\right )_{1\leq i,j\leq n}\!\!=\,
\frac{\theta _{1}^{n-1}(\zeta)\theta _{1}(\zeta
+2\sum_{i=1}^{n}x_i)}{\prod_{i=1}^{n}\theta_{1}(2x_i)}\,
\prod_{i<j}^{n}
\frac{\theta_{1}^{2}(x_i-x_j)}{\theta_{1}^{2}(x_i+x_j)}
$$
which is a particular case of the formula for the elliptic
Cauchy determinant. Let $\Lambda$
be the set $\{1,2,\ldots ,\ell\}$. We use the following
notation. For any subset
$J\subseteq \Lambda$,
$\Lambda \setminus J$ is its complement,
$|J|$ is the number of its elements, and
$\Vert J \Vert =\sum_{m\in J}m$.
Setting $x_n=-2\eta n$, we have:
\begin{eqnarray}
\label{e1}
&&\det \Bigl ( K^{2m}\delta _{mn}+G_{mn}(\zeta )
\Bigr )_{1\leq m,n\leq \ell}
\nonumber\\
&=& \sum _{j=0}^{N}
\frac{\theta _{1}(\zeta -4j\eta)}{\theta _{1}(\zeta)}
K^{2(N-j)}
\sum_{J\subseteq \Lambda , \Vert J \Vert =j} (-1)^{\kappa (J)}
\prod_{k\in J}\prod_{k'\in \Lambda \setminus J}
\frac{\theta_{1}(2(k+k')\eta)}{\theta_{1}(2(k-k')\eta)}\,,
\end{eqnarray}
where
$$
\kappa (J)=|J|\ell +\frac{1}{2}|J|(|J|-1)\,.
$$
Thus, the coefficient $C_{j}^{(\ell)}$ reads
\beq
\label{Cj}
C_{j}^{(\ell)}=
\sum_{J\subseteq \Lambda , \Vert J \Vert =j}
(-1)^{\kappa (J)+j}
\prod_{k\in J}\prod_{k'\in \Lambda \setminus J}
\frac{\theta_{1}(2(k+k')\eta)}{\theta_{1}(2(k-k')\eta)}\,,
\eeq
and the symmetry $j\leftrightarrow N-j$ is evident.
\square

\noindent
Note that the sum in (\ref{Cj}) runs
over partitions of the number $j$ into {\it distinct} parts not
exceeding $\ell$. Examples are given in Appendix B.

The meaning of (\ref{eq4}) is the same as is explained
after equations (\ref{db7}): it defines a covering of
the complex plane invariant under the map
(\ref{db8}). This allows us to define the
corresponding factor-curve which is precisely $\Gamma_e$.
Therefore, $\Gamma_e$ is a ramified covering of the elliptic
curve with the modular parameter $\tau$.

It is easily seen from (\ref{eq4}) that $\Gamma_{e}$
is invariant under the involution
$(\zeta , \,K) \longmapsto (4N\eta -\zeta , \,K^{-1})$.
Eq.\,(\ref{eq1}) allows one to express the function $E$
through $\zeta$, $K$ by the following formula:
\beq
\label{E}
E=
\frac{\theta _{1}(4\ell \eta )}
{\theta _{1}(2\ell \eta )} \,
\frac{\Psi \bigl ( (\ell \!-\!1)\eta \bigr )}
{\Psi \bigl ( \ell \eta \bigr )} =
-\,\frac{s_{\ell -1}}{s_{\ell}}K^{-1}
\frac{[1][2\ell ]}{[\ell ][\ell \!-\!1]}\,.
\eeq
The coefficients $s_{\ell}$, $s_{\ell -1}$ are given by
the corresponding minors of the matrix
$K^{2m}\delta_{mn}+G_{mn}(\zeta )$. It can be shown that
$E$ is invariant under the above involution
of $\Gamma_e$, so this involution coincides with (\ref{db9}).
In terms of the $\Psi$-function, the hyperelliptic
involution takes $\Psi(x)$ to $\Psi (-x)$. Note also that
eq.\,(\ref{eq4})
defines a singular curve. Indeed,
the fixed points of the hyperelliptic involution are
singular points of the curve $\Gamma_e$, i.e.,
both the $\zeta$- and $K$-derivatives
of the left hand side of (\ref{eq4}) at these points equal
to zero. In the neighbourhoods of these points
different sheets of the curve intersect.
The function $E$ takes different values on these sheets
(which are obtained by resolving the indeterminacy in (\ref{E})),
so it resolves the singularities of the curve.

\paragraph{Remark}
The function $\Psi (x)$ is the common
eigenfunction for all the commuting operators
$A_{\lambda}$ (\ref{op1}). Indeed, commutativity of
$L$ and $A_{\lambda}$ implies that
$\tilde \Psi _{\lambda}(x)=A_{\lambda}\Psi (x)$ is
an eigenfunction of $L$ with the same eigenvalue $E$.
By Theorem 3.2, $\tilde \Psi$ is
proportional to $\Psi$: $\tilde \Psi _{\lambda}(x)
=g(\lambda )\Psi (x)$.
From the symmetry (\ref{symm}) and the normalization
condition $A_{\ell \eta} =1$ we have
$$
A_{\lambda}\Psi (x)=\frac{\Psi (\lambda)}{\Psi (\ell \eta )}
\Psi (x)\,.
$$

Let us conclude this section by examining the behaviour
of the spectral curve in the vicinity of its "infinite points",
i.e., the points at which the function $E$ has poles.
From either (\ref{db7}) or (\ref{eq4}), (\ref{E}) we
conclude that there are two such points:
$\infty_{+}=(\zeta \to 0,\, K\to \infty ,\, E\to \infty )$
and $\infty_{-}=
(\zeta \to 4N\eta ,\, K\to 0 ,\, E\to \infty )$.
In the neighbourhood of $\infty _{\pm}$ we have
$E=K^{\pm 1} +o(K^{\pm 1})$, while the leading terms
of the function $K$ are:
$$
\begin{array}{l}
K^2 =- \displaystyle{\frac{1}{\theta_1(\zeta)}}
\displaystyle{\frac{[\ell ][\ell +1]}{[1]}} +O(1)\,,
\;\;\;\;\; \zeta \to 0\,, \\ \\
K^2 =\theta_1(\zeta -4N\eta )
\displaystyle{\frac{[1]}{[\ell ][\ell +1]}} +o(\zeta -4N\eta )\,,
\;\;\; \zeta \to 4N\eta \,.
\end{array}
$$
The Baker-Akhiezer function
$$
\Psi _{BA}(x,P) =\frac{\Psi (x,P)}{\Psi (\ell \eta , P)}\,,
\;\;\;\;\; P=(\zeta , K, E)\in \Gamma \,,
$$
is easily seen to
have the following asymptotics as $P\rightarrow \infty _{\pm}$:
\beq
\Psi _{BA}(x,P) =K^{\frac{x}{\eta}\pm \ell}
\Bigl ( \xi_{0}^{\pm}(x)+O(K^{\mp 1})\Bigr )
\eeq
with some functions $\xi_{0}^{\pm}(x)$, i.e.,
the pole divisor of the Baker-Akhiezer function is
concentrated at $\infty _{\pm}$.

\section{Explicit hyperelliptic realizations}

In this section we obtain the
equation of the spectral curve of the difference Lam\'e
operator $L$, which has the explicit hyperelliptic
form. This equation
contains two variables: $E$ and $z$. The latter is
the eigenvalue of the operator $A_{(\ell +1)\eta}\equiv A$,
which commutes
with $L$, on their common eigenfunction $\Psi$:
$L\Psi =E\Psi$, $A\Psi =z\Psi$.
Recall that the operator $\bar A =A_{-(\ell +1)\eta}$
commutes with both $A$ and $L$.
Let us write out the trivial identity
$A^2 -(A+\bar A)A +A\bar A =0$ and act by both sides
on the $\Psi$. Taking into account
(\ref{TAA}), we get
$z^2 -([2\ell ]!)^{-1}T_{2\ell +1}(0,E)+D_{2\ell}(E)=0$, or
\beq
\label{H1}
z+\frac{D_{2\ell}(E)}{z}=
([2\ell ]!)^{-1}T_{2\ell +1}(0,E)\,,
\eeq
where
$D_{2\ell}(E)$ is a
polynomial of $E$ of degree $2\ell$ and
$T_{2\ell +1}(0,E)=\lim _{u\to 0}T_{2\ell +1}(u,E)$
is the polynomial of $E$ of degree $2\ell +1$.
They enjoy the properties $D_{2\ell}(-E)=D_{2\ell}(E)$,
$T_{2\ell +1}(0,-E)=-T_{2\ell +1}(0,E)$.
Recall that
\beq
\label{LA}
L=c_{-}(x)e^{\eta \p _{x}}
+c_{+}(x)e^{-\eta \p _{x}}\,,
\;\;\;\;\;\;\;
A=\sum_{k=0}^{\ell}a_{2k+1}(x)e^{(2k+1)\eta \p _{x}}\,,
\eeq
where the coefficient functions are:
\beq
\label{H2}
\begin{array}{l}
c_{\pm}(x)=\displaystyle{
\frac{\theta_1(2x \pm \ell \eta)}{\theta_1(2x)}}\,,
\\ \\
a_{2k+1}(x)=(-1)^k\displaystyle{\frac{
\left [ \begin{array}{c}2\ell +1 \\ \ell -k
\end{array}\right ]}{\left [ \begin{array}{c}2\ell \\
\ell \end{array}\right ]}}
\prod _{j=0}^{\ell -k -1}
\displaystyle{\frac{\theta_1(2x+2(\ell
-j)\eta )}{\theta_1(2x+2(\ell \!+\!k\!-\!j\!+\!1)\eta )}}
\prod _{j=0}^{k-1}\displaystyle{\frac{\theta_1(2x-2(\ell
-j)\eta )}{\theta_1(2x+2(k\!+\!j\!+\!1)\eta )}}
\end{array}
\eeq
(see (\ref{op1}), (\ref{op2})).

Introducing $w=2z-([2\ell ]!)^{-1}T_{2\ell +1}(0,E)$,
we rewrite (\ref{H1}) in the customary hyperelliptic form:
\beq
\label{H1a}
w^2 =([2\ell ]!)^{-2}(T_{2\ell +1}(0,E))^2-4D_{2\ell}(E)
=\left [ \begin{array}{c} 2\ell \\ \ell \end{array}\right ]^{-2}
P_{2\ell +1}(E^2)\,,
\eeq
where $P_{2\ell +1}(E^2)=\prod _{i=1}^{2\ell +1}(E^2 -E_{i}^{2})$.
Note that $w$ is the eigenvalue of the operator
$W=A-\bar A$.

To find out the explicit form of $D_{2\ell}(E)$,
we make use of the following simple argument:
given a commuting pair $[L,A]=0$ of
operators of finite order,
the spectral problem for $L$ is reduced to an
eigenvalue problem for a finite matrix with elements depending on the
eigenvalue of $A$ (the "spectral parameter").
Specifically, let $\Psi (x)$ be their common eigenfunction.
We set
\beq
\label{H3}
\Psi _j =\Psi_j(x)=\Psi
(x+j\eta)\,.
\eeq
The equation $A\Psi =z\Psi$ allows us to express $\Psi_0$
and $\Psi_{2\ell+2}$ through
$\Psi_1, \Psi_2, \ldots , \Psi _{2\ell +1}$:
\beq
\label{H4}
\begin{array}{l}
\Psi_0 = z^{-1}
\displaystyle{\sum_{k=0}^{\ell}}a_{2k+1}(x)\Psi_{2k+1}\,, \\ \\
\Psi_{2\ell +2} = za^{-1}_{2\ell +1}(x+\eta)\Psi_1 -
\displaystyle{\sum_{k=1}^{\ell}}\displaystyle{
\frac{a_{2k-1}(x+\eta )}{a_{2\ell +1}(x+\eta )}}
\Psi_{2k}\,.
\end{array}
\eeq
Now the spectral problem $L\Psi =E\Psi$ can be rewritten as a
homogeneous
linear system for $\Psi_1, \ldots , \Psi _{2\ell +1}$.
Equating its determinant to zero, we get a relation
between $z$ and $E$, which is the equation of our spectral
curve. Its independence of the value of $x$ follows,
eventually, from commutativity of $L$ and $A$.

To be more precise, introduce the vector-function
$\vec \Psi (x)$ with components $\Psi_j(x)$
(see (\ref{H3})), $j=1,\ldots , 2\ell +1$.
The eigenvalue equation for $A$ can be rewritten in the form
\beq
\label{H5}
\vec \Psi (x-\eta )={\bf A}(x,z)\vec \Psi (x)\,,
\eeq
where the matrix ${\bf A}(x,z)$ reads
\beq
\label{H6}
{\bf A}(x,z)=\left (
\begin{array}{ccccccccccc}
z^{-1}a_1(x)&&0&&z^{-1}a_3(x)&&\ldots && 0 &&
z^{-1}a_{2\ell \!+\!1}(x) \\ &&&&&&&&&& \\
1 && 0 && 0 && \ldots && 0 && 0 \\ &&&&&&&&&& \\
0 && 1 && 0 && \ldots && 0 && 0 \\ &&&&&&&&&& \\
\ldots && \ldots && \ldots && \ldots &&
\ldots  && \ldots \\ &&&&&&&&&& \\
\ldots && \ldots && \ldots && \ldots &&
\ldots  && \ldots \\ &&&&&&&&&& \\
0 && 0 && 0 && \ldots && 1 && 0
\end{array} \right )\,.
\eeq
The scalar equation $L\Psi =E\Psi$ is then rewritten
in the matrix form ${\bf L}(x,z)\vec \Psi (x) =E\vec \Psi (x)$,
where
\beq
\label{H7}
{\bf L}(x,z)={\bf C}_{+}(x){\bf A}(x,z)
+{\bf C}_{-}(x){\bf A}^{-1}(x+\eta ,z)\,,
\eeq
$$
{\bf C}_{\pm}=\mbox{diag}\,\Bigl \{
c_{\pm}(x+\eta ), \, c_{\pm}(x+2\eta ), \ldots ,
c_{\pm}(x+(2\ell +1)\eta ) \Bigr \}\,.
$$
So, we have assigned to the scalar operators $A$ and $L$
the $(2\ell +1)\times (2\ell +1)$-matrices ${\bf A}$ and
${\bf L}$ respectively.
\begin{lem}
The commutativity condition $[L,A]=0$ is equivalent
to the Lax equation
\beq
\label{H8}
{\bf L}(x-\eta , z){\bf A}(x,z)={\bf A}(x,z){\bf L}(x, z)
\eeq
for arbitrary values of $z$.
\end{lem}
{\it Proof.} This equality holds identically for the rows
from the second to the last one. The first row gives the set
of relations for coefficients of the $L,A$, which are
equivalent to their commutativity.
\square

\noindent
The system
${\bf L}(x,z)\vec \Psi (x) =E\vec \Psi (x)$
has nontrivial solutions if and only if $E$ and $z$ are
connected by the equation of the spectral curve:
\beq
\label{H9}
\det \Bigl ( {\bf L}_{ij}(x,z)-
E\delta_{ij}\Bigr )_{1\leq i,j\leq 2\ell +1}\!=0\,.
\eeq
It is easy to see that this determinant does not depend on $x$.
Indeed, denote the left hand side of (\ref{H9}) by
$f(x)$. It follows from (\ref{H7}) that
${\bf L}_{mn}(x+\frac{1}{2},z)={\bf L}_{mn}(x,z)$,
${\bf L}_{mn}(x+\frac{1}{2}\tau ,z)=
e^{-4\pi i \ell (m-n)\eta }{\bf L}_{mn}(x,z)$,
whence $f(x)$ is a double-periodic function of the variable
$y=2x$ with periods $1$ and $\tau$ and with finite number of
possible poles. At the same time, (\ref{H8}) implies
that $f(x)$ has one and the same value at
infinite number of points $x+m\eta$, $m\in {\bf Z}$.
Thus, $f(x) =\mbox{const}$.

Extracting the $z$-dependence of the
determinant in (\ref{H9}), we can write
$$
\det \Bigl ( {\bf L}_{ij}(x,z)-
E\delta_{ij}\Bigr )= z(-1)^{\ell}\left [
\begin{array}{c} 2\ell \\ \ell \end{array}\right ]
+ F + z^{-1}G\,,
$$
where $F$ and $G$ are polynomials of $E$. According to the above
argument, they do not depend on $x$. It is convenient to
evaluate $F$ at $x\to -(\ell +1)\eta$. This should be done
with some care because some matrix elements are singular at this
point. Using (\ref{T2ell1}), we find:
$$
\begin{array}{lll} F&=&
\left ( \displaystyle{ \prod _{k=1, \neq \ell+1}^{2\ell +1}}
\theta_1(2x+2k\eta )\right )^{-1}T_{2\ell +1}
(2x\!+\!2(\ell \!+\!1)\eta , -E)  \\ && \\
&+& \mbox{two "unwanted" determinants}\,.
\end{array}
$$
At $x\to -(\ell +1)\eta$ the first term yields
$(-1)^{\ell +1} ([\ell ]!)^{-2}T_{2\ell +1}(0,E)$, and
each of the two "unwanted" terms tends to zero.
The simplest determinant representation for
$G$ is obtained at $x=\ell \eta$ when almost all
elements of the first row are equal to zero.
Skipping the details,
we present the result.

The equation of the curve has the form (\ref{H1}), where
\beq
\label{H10}
D_{2\ell}(E)=(-1)^{\ell}\frac{[2\ell +1]}{[\ell +1]}
\left [ \begin{array}{c}2\ell \\ \ell \end{array}\right ]^{-1}
\det ({\bf D}_{ij})_{1\leq i,j\leq 2\ell}\,.
\eeq
The $(2\ell \times 2\ell )$-matrix ${\bf D}_{ij}$ reads
\beq
\label{H11}
{\bf D}_{ij}=\left (
\begin{array}{ccccccccccccccc}
-E&&
\frac{[2]}{[\ell \!+\!2]}&& 0 && 0 && \ldots && 0 && 0 && 0\\
&&&&&&&&&&&&&&\\
\frac{[2\ell \!+\!3]}{[\ell \!+\!3]}&& -E &&
\frac{[3]}{[\ell \!+\!3]}&& 0 && \ldots && 0 && 0 && 0\\
&&&&&&&&&&&&&&\\
0 &&
\frac{[2\ell \!+\!4]}{[\ell \!+\!4]}&& -E &&
\frac{[4]}{[\ell \!+\!4]}&& \ldots && 0 && 0 && 0\\
&&&&&&&&&&&&&&\\
\ldots && \ldots && \ldots && \ldots &&
\ldots && \ldots && \ldots && \ldots \\
&&&&&&&&&&&&&&\\
\ldots && \ldots && \ldots && \ldots &&
\ldots && \ldots && \ldots && \ldots \\
&&&&&&&&&&&&&&\\
0 && 0 && 0 && 0 && \ldots &&
\frac{[4\ell ]}{[3\ell ]}&& -E &&
\frac{[2\ell ]}{[3\ell ]}\\
&&&&&&&&&&&&&&\\
d_2 && 0 && d_4 && 0 && \ldots && 0 &&
d_{2\ell}\!+\!
\frac{[4\ell \!+\!1]}{[3\ell \!+\!1]}&& -E
\end{array}
\right ),
\eeq
where the coefficients $d_{2k}$ entering the last line are:\
$$
d_{2k}=(-1)^{\ell -k}\frac{[\ell +2k]}{[k]}
\left [ \begin{array}{c}3\ell \\ \ell \end{array}\right ]
\left [ \begin{array}{c}2\ell +1 \\ \ell +k\end{array}\right ]
\left [ \begin{array}{c}2\ell \!+\!k\!+\!1
\\ k\end{array}\right ]^{-1},
\;\;\;\;k=1,\ldots , \ell\,.
$$
It is not difficult to see that the element
${\bf D}_{2\ell , 2\ell -1}$ is equal to
$$
d_{2\ell}+
\frac{[4\ell \!+\!1]}{[3\ell \!+\!1]}=
\frac{[\ell +1][2\ell ]^2}{[1][\ell ]^2}\,.
$$
The polynomials $D_{2\ell}(E)$ for small values
of $\ell$ are written out in Appendix B.

The size of the determinant in the equation defining the
curve can be reduced for the price of more complicated
matrix elements. Let us mention two examples.

The first one is 2$\times$2 determinant
representation, which is in certain sense "dual" to (\ref{H8}).
The duality means that the roles
of ${\bf A}$- and ${\bf L}$-matrices are interchanged:
according to the same scheme,
the spectral problem for $A$ is reduced to an eigenvalue
problem for a 2$\times$2-matrix with elements depending
on the eigenvalue $E$ of $L$. Specifically, set
\beq
\label{H12}
{\cal L}(x,E)=\left (
\begin{array}{ccc}
0 && 1 \\ && \\
-\displaystyle{\frac{c_{+}(x+\eta )}{c_{-}(x+\eta )}} &&
Ec_{-}^{-1}(x+\eta )
\end{array}
\right )
\eeq
and rewrite the equation $L\Psi =E\Psi$ in the form
$\left ( \begin{array}{c}\Psi (x+\eta ) \\ \Psi (x+2\eta )
\end{array} \right ) ={\cal L}(x,E)
\left ( \begin{array}{c}\Psi (x) \\ \Psi (x+\eta )
\end{array} \right )$.
Introduce the "monodromy matrix"
$$
{\cal A}(x,E)=\sum _{k=0}^{\ell}\left (
\begin{array}{cc} a_{2k+1}(x)& 0 \\
0 & a_{2k+1}(x+\eta ) \end{array} \right )
\prod _{2k\geq j \geq 0}^{\longleftarrow} {\cal L}(x+j\eta )\,,
$$
where the arrow indicates the ordered product of matrices,
and $a_{2k+1}(x)$ are the same as in (\ref{LA}).
Then the spectral problem $A\Psi =z\Psi$ is rewritten
as ${\cal A}(x,E)
\left ( \begin{array}{c}\Psi (x) \\ \Psi (x+\eta )
\end{array} \right ) =z
\left ( \begin{array}{c}\Psi (x) \\ \Psi (x+\eta )
\end{array} \right )$, and the Lax equation
\beq
\label{H13}
{\cal A}(x+\eta , E){\cal L}(x,E)
={\cal L}(x, E){\cal A}(x,E)
\eeq
holds true provided the operators $L$ and $A$ commute.
Note that now ${\cal A}(x,E)$ plays the role
of the Lax matrix (cf. (\ref{H8})). The equation of the
spectral curve reads
\beq
\label{H14}
\det ({\cal A}(x,E) -z)=0\,.
\eeq
Due to the Lax equation (\ref{H13}) it does not
depend on $x$.

The second one is $\ell \times \ell$ determinant representation
for the factor-curve obtained by factorization over the reflection
$E\rightarrow -E$. In the operator language,
this amounts to finding a polynomial relation between the
commuting operators $A_{(\ell -2)\eta}$ and
$A_{(\ell +2)\eta}$. As it follows from (\ref{bax1}), their
eigenvalues ($\varepsilon$ and
$\xi$, respectively) are connected with $E,z$ by the formulas
$$
\varepsilon =\frac{[\ell -1][\ell ]}{[2\ell ][2\ell -1]}
\left ( E^2 +
\frac{[1][2\ell ]}{[\ell -1][\ell ]}\right ),
\;\;\;\;\;
\xi =\frac{[\ell +1]}{[1]}
\left ( Ez -
\frac{[2\ell +1]}{[\ell +1]}\right ).
$$
The equation of the curve can be
derived in the same way as (\ref{H9}).

\section{Remarks on representations
of the Sklyanin algebra}

In this section we make a few remarks relating the
above material to representation theory of the Sklyanin
algebra.

Let ${\cal S}_{a}^{(\ell )}$ denote the generators of the
Sklyanin algebra realized by difference operators as in
(\ref{Sa}). To make connections between the commuting family
$A_{\lambda}$ (\ref{op1}) and representations of the Sklyanin
algebra explicit,
we begin with a simple reformulation of the
Novikov equation $[L,\, V]=0$
for coefficients of an operator $V$ commuting
with $L={\cal S}_{0}^{(\ell )}$.
Suppose $v(x , y)$ is any solution to the equation
\beq
\label{01}
\nabla v(x , y) =0\,,
\eeq
where $\nabla ={\cal S}_{0}^{(\ell )}(x,\p_x)-
{\cal S}_{0}^{(\ell )}(y,\p_y)$. (The operator
${\cal S}_{0}^{(\ell )}(x,\p_x)$ acts to the variable $x$,
etc.)
Then the operator
\beq
V=\sum _{}V_j(x)
e^{j\eta \p _x}
\label{sum}
\eeq
with the coefficients
$$
V_j(x)=v(x , x +j\eta )\left (
\prod _{k=1}^{\ell}\theta _{1}(2x+2(j-k)\eta )
\theta _{1}(2x+2(j+k)\eta )\right )^{-1}
$$
commutes with ${\cal S}_0^{(\ell )}$.
In general, the operators $V$ given by this construction
are of infinite order, i.e. the sum in (\ref{sum}) is infinite.
For the family of commuting operators
$A_{\lambda}$ (\ref{op1}) the sum is {\it finite}.
This corresponds to some very particular solutions to
eq.\,(\ref{01}).

Recall that the Sklyanin algebra can be realized \cite{krizab}
by certain difference operators in two variables acting on the
invariant subspace of solutions to eq.\,(\ref{01}).
Taking into account the
isomorphism between difference
operators commuting with the operator
${\cal S}_{0}^{(\ell )}$ and meromorphic functions on
its spectral curve $\Gamma$ with poles only
at $\infty _{\pm}$, we conclude that the Sklyanin algebra
acts in the space of such functions on $\Gamma$.
It would be very interesting to find the explicit form
of this action.

In the rest of this section we present some results
on the role of the commuting operators $A_{\lambda}$
in representations of the Sklyanin algebra.
First, we show that they have the same invariant
subspace as the generators of the Sklyanin algebra (\ref{Sa}).
Second, the operator $W=A\!-\!\bar A$ is shown to "intertwine"
representations of spins $\ell$ and $-\ell -1$.

\begin{prop}
Let $\ell$ be a positive integer. Then
the operators $A_{\lambda}$ (\ref{op1})
preserve the space $\Theta_{4\ell}^{+}$.
\end{prop}
{\it Sketch of proof.}
Let $F(x)\in \Theta_{4\ell}^{+}$. Then the monodromy properties
of $\tilde F(x) =
(A_{\lambda}F)(x)$ are the same as thous of
$\theta$-functions of order $4\ell$ (this is easily seen from
(\ref{op2})). Next, a
further inspection of (\ref{op1}), (\ref{op2}) shows
that the condition $F(x)=F(-x)$
is enough for cancellation of all poles
of $\tilde F(x)$. Therefore, $\tilde F(x)\in \Theta _{4\ell}$.
It remains to prove that $\tilde F(x)$ actually belongs to
$\Theta^{+}_{4\ell}$. Set $f(x)=\tilde F(x)\!-\!\tilde F(-x)$,
then $f(x)=A_x (\lambda , \p _{\lambda})F(\lambda)-
A_{-x}(\lambda , \p _{\lambda})F(\lambda)$, where the
$x\leftrightarrow \lambda$ symmetry (\ref{symm}) is used.
Since $f(x)\in \Theta _{4\ell}$, it is enough to verify
the equality $f(x)=0$ in $4\ell$ points $x=m\eta +\omega _{a}$,
$m=1, \ldots , \ell$, where $\omega_a$ are the half-periods.
This is easy to do if to recall the operator identity
$A_{m\eta}=A_{-m\eta}$ for $m=1, \ldots , \ell$.
\square
\begin{cor}
For any $\lambda \in {\bf C}$ the operator
$A_{\lambda}\!-\!A_{-\lambda}$
annihilates the space $\Theta_{4\ell}^{+}$.
\end{cor}
Indeed, for $F(x) \in \Theta _{4\ell}^{+}$ the function
$(A_{\lambda}-A_{-\lambda})F(x)$ is simultaneously odd and even.
\paragraph{Remark}
Recall the involution (\ref{db81})
that changes the sign of $E$. This involution takes
the eigenfunction $\Psi (x)$ (see \ref{psi1}) to
$e^{\frac{i\pi x}{\eta}} \Psi (x)$. Clearly,
the operators $A_{\lambda}$ preserve
the space $e^{\frac{i\pi x}{\eta}} \Theta _{4\ell}^{+}$ as well,
so $A_{\lambda}\!-\!A_{-\lambda}$ annihilates the space
$\Theta _{4\ell}^{+}\oplus e^{\frac{i\pi x}{\eta}}
\Theta _{4\ell}^{+}$ spanned by eigenfunctions of the
difference Lam\'e operator at the edges of bands (cf. \cite{kz}).

To formulate the next proposition, it is convenient to
slightly modify the operator $W=A-\bar A$. Let us
introduce the operator
\beq
\check W=(-1)^{\ell}\left [ \begin{array}{c} 2\ell \\ \ell
\end{array} \right ]\varphi ^{-1}_{\ell}(x) W\,,
\label{S1}
\eeq
where
$\varphi _{\ell}(x)=\prod _{j=0}^{2\ell }
\theta_1(2x+2(j-\ell )\eta )$.
The explicit formula for $\check W$ can be written
in the form that has sense
not only for $\ell \in {\bf Z}_{+}$ but also for
$\ell \in {\bf Z}_{+}+\frac{1}{2}$:
\beq
\label{S2}
\check W = \sum_{k=0}^{2\ell +1} (-1)^k
\left [ \begin{array}{c} 2\ell +1 \\ k \end{array} \right ]
\frac{ \theta _{1}(2x+
2(2\ell -2k +1)\eta )}{\prod_{j=0}^{2\ell -k+1}
\theta_1(2x+2j\eta )
\prod_{j'=1}^{k}\theta_1(2x-2j'\eta )}
e^{(2\ell -2k +1)\eta \p_x }\,.
\eeq
Here and below
the dependence of the $\check W$ on $\ell$ is not
indicated explicitly. The following proposition is proved
by a straightforward verification using some identities
for the $\theta$-functions.
\begin{prop}
For $\ell \in \frac{1}{2}{\bf Z}_{+}$,
the operator $\check W$
"intertwines" representations
of spin $\ell$ and of spin $-(\ell +1)$:
\beq
\label{S3}
{\cal S}_{a}^{(-\ell -1)}\check W =\check W
{\cal S}_{a}^{(\ell )}\,,
\;\;\;\;\;\; a=0,\ldots , 3\,.
\eeq
\end{prop}
The same intertwining relation can be written
for the quantum $\mbox{{\sf L}}$-operator (\ref{L}):
$\mbox{{\sf L}}^{(-\ell -1)}\check W =
\check W \mbox{{\sf L}}^{(\ell )}$.
\paragraph{Remark}
In case of the algebra $sl_2$ the intertwining
operator between representations of spins $\ell$ and
$-\ell -1$ (realized by differential operators in $x$) is
$(d/dx)^{2\ell +1}$. It annihilates the linear space
of polynomials of degree $\leq 2\ell$.

\noindent
Note that the operator $\check W$ is not invertible.
By Corollary 5.1, for $\ell \in {\bf Z}_{+}$
$\check W$ annihilates the space
$\Theta _{4\ell}^{+}\oplus e^{\frac{i\pi x}{\eta}}
\Theta _{4\ell}^{+}$ (see the remark after Corollary 5.1).
As is mentioned above, this is precisely the space
spanned by eigenfunctions of the
difference Lam\'e operator at the edges of bands. So, in this way
we obtain another proof of the result of \cite{kz}: the eigenfunctions
of $L$ at the edges of bands span a $(4\ell +2)$-dimensional
functional subspace, which is invariant for all Sklyanin's
operators ${\cal S}_a$. The corresponding
$(4\ell +2)$-dimensional representation of the Sklyanin
algebra is the direct sum of two equivalent
$(2\ell +1)$-dimensional irreducible representations.

\section*{Acknowledgements}
The author thanks I.Krichever, A.Mironov and T.Takebe
for useful discussions.
This work was supported in part by RFBR grant 98-01-00344
and by grant 96-15-96455 for support of scientific schools.

\section*{Appendix A. Theta-functions}
\def\theequation{A\arabic{equation}}
\setcounter{equation}{0}

We use the following definition of the
Jacobi $\theta$-functions:
\beq
\begin{array}{l}
\theta _1(x|\tau)=\displaystyle{\sum _{k\in {\bf Z}}}
\exp \left (
\pi i \tau (k+\frac{1}{2})^2 +2\pi i
(x+\frac{1}{2})(k+\frac{1}{2})\right ),
\\ \\
\theta _2(x|\tau)=\displaystyle{\sum _{k\in {\bf Z}}}
\exp \left (
\pi i \tau (k+\frac{1}{2})^2 +2\pi i
x(k+\frac{1}{2})\right ),
\\ \\
\theta _3(x|\tau)=\displaystyle{\sum _{k\in {\bf Z}}}
\exp \left (
\pi i \tau k^2 +2\pi i
xk \right ),
\\ \\
\theta _4(x|\tau)=\displaystyle{\sum _{k\in {\bf Z}}}
\exp \left (
\pi i \tau k^2 +2\pi i
(x+\frac{1}{2})k\right ).
\end{array}
\label{theta}
\eeq
Throughout the paper we write
$\theta _a(x|\tau)=\theta _a(x)$. The frequently used
transformation properties for shifts by (half) periods are:
\beq
\label{periods}
\theta_a (x\pm 1)=(-1)^{\delta _{a,1}+\delta _{a,2}}
\theta_a (x)\,,
\;\;\;\;\;
\theta_a (x\pm \tau )=(-1)^{\delta _{a,1}+\delta _{a,4}}
e^{-\pi i \tau \mp 2\pi i x}
\theta_a (x)\,,
\eeq
\beq
\label{half}
\begin{array}{l}
\theta_1 (x\pm \frac{1}{2})=\pm
\theta_2 (x)\,,\\ \\
\theta_1 (x\pm \frac{\tau}{2})=\pm i
e^{-\frac{1}{4}\pi i \tau \mp \pi i x}
\theta_4 (x)\,, \\ \\
\theta_1 (x\pm \frac{1+\tau}{2})=\pm
e^{-\frac{1}{4}\pi i \tau \mp \pi i x}
\theta_3 (x)\,,
\end{array}
\eeq

By $\Theta_n$ we denote the space of $\theta$-functions
of order $n$, i.e., entire functions
$F(x)$, $x\in {\bf C}$, such that
\beq
F(x+1)=F(x)\,,
\;\;\;\;\;\;
F(x+\tau)=(-1)^n e^{-\pi i n\tau -2\pi i nx}F(x)\,.
\label{8}
\eeq
It is easy to see that $\mbox{dim} \,\Theta_n =n$.
Let $F(x)\in \Theta_n$, then $F(x)$ has a multiplicative
representation of the form
$$
F(x)=c\prod _{i=1}^{n}\theta_1(x-x_i)\,,
\;\;\;\;\;\;\;\sum _{i=1}^{n}x_i =0\,,
$$
where $c$ is a constant. Imposing, in addition to (\ref{8}),
the condition $F(-x)=F(x)$, we define the space
$\Theta_{n}^{+}\subset \Theta_{n}$ of {\it even}
$\theta$-functions of order $n$, which
plays the important role in representations
of the Sklyanin algebra. If $n$ is an even number,
then $\mbox{dim}\, \Theta_{n}^{+} =\frac{1}{2}n +1$.

\section*{Appendix B }
\def\theequation{A\arabic{equation}}
\setcounter{equation}{0}

Here we explicitly write out equation (\ref{eq4}) for
small values of $\ell$. The "eliptic number" notation
$[n]\equiv \theta_1(2n\eta )$ is used.
$$
\begin{array}{lll} \ell =1:&&
\!\theta_1 (\zeta)K^2 -\theta_1 (\zeta \!-\!4\eta)=0\,,\\
&&\\
\ell =2: &&
\!\theta_1 (\zeta)K^6 -
\displaystyle{\frac{[3]}{[1]}}
\theta_1(\zeta \!-\!4\eta)K^4
+\displaystyle{\frac{[3]}{[1]}}
\theta_1(\zeta \!-\!8\eta)K^2
-\theta_1 (\zeta \!-\!12\eta)=0\,,\\
&&\\
\ell =3: &&
\!\!\theta_1 (\zeta)K^{12} -
\displaystyle{\frac{[3][4]}{[1][2]}}
\theta_1(\zeta \!-\!4\eta)K^{10}
+\displaystyle{\frac{[3][5]}{[1]^2}}
\theta_1(\zeta \!-\!8\eta)K^{8}-
2\,\displaystyle{\frac{[4][5]}{[1][2]}}
\theta_1(\zeta \!-\!12\eta)K^{6}+\\
&&\\
\phantom{\ell =3:}&&
+\displaystyle{\frac{[3][5]}{[1]^2}}
\theta_1(\zeta \!-\!16\eta)K^{4}
-\displaystyle{\frac{[3][4]}{[1][2]}}
\theta_1(\zeta \!-\!20\eta)K^{2}
+\theta_1 (\zeta -24\eta)=0\,,
\end{array}
$$
$$
\begin{array}{lll}
\ell =4: &&
\!\!\!\!\theta_1 (\zeta)K^{20} \!-\!
\displaystyle{\frac{[4][5]}{[1][2]}}
\theta_1(\zeta \!-\!4\eta)K^{18}
\!+\!\displaystyle{\frac{[3][5][6]}{[1]^2[2]}}
\theta_1(\zeta \!-\!8\eta)K^{16}-\\
&&\\
\phantom{\ell =4:}&&
\!\!-\left (
\displaystyle{\frac{[4][5][7]}{[1]^2[2]}}
\!+\!
\displaystyle{\frac{[4][5]^2[6]}{[1][2]^2[3]}}\right )
\theta_1(\zeta \!-\!12\eta)K^{14}
\!+\!\left (
\displaystyle{\frac{[5][6][7]}{[1][2][3]}}
\!+\!
\displaystyle{\frac{[5]^2[7]}{[1]^3}}\right )
\theta_1(\zeta \!-\!16\eta)K^{12}-\\
&&\\
\phantom{\ell =4:}&&
-2\,
\displaystyle{\frac{[3][4][6][7]}{[1]^2[2]^2}}
\theta_1(\zeta \!-\!20\eta)K^{10}+\\
&&\\
\phantom{\ell =4:}&&
+\left (
\displaystyle{\frac{[5][6][7]}{[1][2][3]}}
+
\displaystyle{\frac{[5]^2[7]}{[1]^3}}\right )
\theta_1(\zeta \!-\!24\eta)K^{8}
-\left (
\displaystyle{\frac{[4][5][7]}{[1]^2[2]}}
+
\displaystyle{\frac{[4][5]^2[6]}{[1][2]^2[3]}}\right )
\theta_1(\zeta \!-\!28\eta)K^{6}+\\
&&\\
\phantom{\ell =4:} &&
+\displaystyle{\frac{[3][5][6]}{[1]^2[2]}}
\theta_1(\zeta \!-\!32\eta)K^{4}-
\displaystyle{\frac{[4][5]}{[1][2]}}
\theta_1(\zeta \!-\!36\eta)K^{2}+
\theta_1(\zeta \!-\!40\eta)=0\,.
\end{array}
$$

Let us present the explicit form of the polynomials
$T_{2\ell +1}(0,E)$ and $D_{2\ell}(E)$ (see (\ref{H1})) for
$\ell =1,2$.
$$
\begin{array}{lll}
T_3(0,E)&=&-[1]^2\left \{ E^3 -
\left ( \displaystyle{\frac{[1][4]}{[2][3]}}+
\displaystyle{\frac{ [2]^4}{[1]^3[3]}} \right )E\right \}\,,\\
&& \\
T_5(0,E)&=&[1]^2 [2]^2 \left \{ E^5 +
\left ( 3\displaystyle{\frac{[4]}{[2]}}
-\displaystyle{\frac{[2]^4}{[1]^4}}\right )E^3
+ \left (
\displaystyle{\frac{[4]^2}{[2]^2}}+
\displaystyle{\frac{[3]^3}{[1]^3}}-
\displaystyle{\frac{[1] [6]}{[2][3]}}\right )E
\right \}\,,
\end{array}
$$

$$
\begin{array}{lll}
D_2(E)&=&-\displaystyle{\frac{
[1]}{[2]}}\left (
\displaystyle{\frac{ [3]}{[2]}}E^2 -
\displaystyle{\frac{ [2]^3}{[1]^3}}
\right )\,,\\
&&\\
D_4(E)&=&\displaystyle{\frac{[1]}{[3]}}\left \{
\displaystyle{\frac{ [2][5]}{[3][4]}}E^4 -
\left (
\displaystyle{\frac{ [2]^2[7]}{[3][4]^2}}+
\displaystyle{\frac{ [4]^2[5]}{[1][2][6]}}+
\displaystyle{\frac{ [2][8]}{[4][6]}}
\right ) E^2 +
\displaystyle{\frac{ [4][5]}{[1][6]}}
\left (
\displaystyle{\frac{ [7]}{[5]}}
+\displaystyle{\frac{ [5]}{[1]}}
\right ) \right \}\,.
\end{array}
$$

\end{document}